\documentclass[graybox]{svmult}
\usepackage{mathptmx}       
\usepackage{helvet}         
\usepackage{courier}        
\usepackage{type1cm}        
\usepackage{makeidx}         
\usepackage{graphicx}        
\usepackage{multicol}        
\usepackage[bottom]{footmisc}
\usepackage{url}
\usepackage{amssymb}
\usepackage{amsmath}

\def\Prob{\operatorname{Prob}}

\begin{document}
\title*{Piecewise deterministic Markov processes in biological models}
\author{Ryszard Rudnicki and Marta Tyran-Kami\'nska}
\institute{Ryszard Rudnicki \at Institute of Mathematics,
Polish Academy of Sciences, Bankowa 14, 40-007 Katowice, Poland,
\email{rudnicki@us.edu.pl}
\and
Marta Tyran-Kami\'nska \at
Institute of Mathematics,
University of Silesia, Bankowa 14, 40-007 Kato\-wi\-ce, Poland,
\email{mtyran@us.edu.pl}
\and This research was partially supported
by the State Committee for Scientific Research (Poland) Grant No. ~N~N201~608240. The first author is
a supervisor  in the International Ph.D. Projects Programme of Foundation for Polish Science operated within the
Innovative Economy Operational Programme 2007-2013 (Ph.D. Programme: Mathematical Methods in Natural Sciences).
}
\maketitle

\abstract*{We present a short introduction into the framework of piecewise
deterministic Markov processes. We illustrate the abstract mathematical setting with a series of examples
related to dispersal of biological systems, cell cycle models, gene expression, physiologically structured
populations, as well as neural activity. General results concerning asymptotic properties of stochastic
semigroups induced by such Markov processes are applied to specific examples.}

\abstract{We present a short introduction into the framework of piecewise
deterministic Markov processes. We illustrate the abstract mathematical setting with a series of examples
related to dispersal of biological systems, cell cycle models, gene expression, physiologically structured
populations, as well as neural activity. General results concerning asymptotic properties of stochastic
semigroups induced by such Markov processes are applied to specific examples.}


\section{Introduction}
\label{sec:intro} The aim of this chapter is to give a short mathematical introduction to piecewise
deterministic Markov processes (PDMPs) including some results concerning their asymptotic behavior and providing
biological models where they appear. According to a non-rigorous definition by Davis \cite{davis84}, the class
of  \textit{piecewise deterministic Markov processes} is a general family of stochastic models covering
virtually
 all non-diffusion applications.  A more formal definition is the following:
a continuous time Markov process $X(t)$, $t\ge 0$, is a PDMP if  there is an increasing sequence of random times
$(t_n)$, called jump times, such that the sample paths  of $X(t)$ are defined in a deterministic way in each
interval $(t_n,t_{n+1})$. We consider two types of behavior  of the process at jump times: the process can jump
to a new point  or can change the dynamics which defines its trajectories. PDMPs is a large family of different
stochastic processes which includes discrete time Markov processes, continuous time Markov chains, deterministic
processes with jumps, processes with switching dynamics and some point processes. Although the discrete time
Markov processes play important role in applications we will not investigate them here because their theory
differs from that of continuous time PDMPs and their applications  are sufficiently known~\cite{Allen}.

The outline of this chapter is as follows.  In Section 2 we present a number of simple biological models to
illustrate possible applications of such processes.  In Section 3 we collect relevant definitions and examples
of stochastic semigroups. In Section 4 we recall two general results concerning the long-time behavior
(asymptotic stability and sweeping) of stochastic semigroups and we show how they can be applied in the context
of PDMPs with switching dynamics. Examples of applications of these results to concrete biological models are
also provided. The chapter concludes with a short summary and discussion.

\section{Examples}
\label{sec:examples}
\subsection{Pure jump-type and velocity jump Markov processes}
\label{ss:pj-vj} The simplest examples of PDMPs are continuous time Markov chains. Their theory is well known,
so we only mention here that they have a lot of biological applications such as  birth-death processes,
epidemic models (see \cite{Allen})   and, more recently, models of genome evolution  (see e.g. \cite{RTW,RT}).
Continuous time Markov chains belong to a slightly larger class of the so-called pure jump-type Markov
processes. A \textit{pure jump-type Markov process} is a Markov process which remains constant between jumps.
For example, the process used in a simple description of the grasshopper  and kangaroo movement \cite{ODA} is an
example of a pure jump-type Markov process, which is not a Markov chain.  A grasshopper jumps at random times
$t_n$ from a point $x$ to the point $x+Y_n$. We assume that jump times are the same as for a Poisson process
$N(t)$  with intensity $\lambda>0$, i.e., $N(t_n)=n$, and that  $(Y_n)$ is a sequence of independent and
identically distributed (i.i.d.) random vectors. Then the position $X(t)$ of the grasshopper at time $t$ is
given by
\begin{equation}
\label{1-gm1}
X(t)=X(0)+\sum_{n=1}^{N(t)} Y_n.
\end{equation}
The process as in  (\ref{1-gm1}) is called a \textit{compound Poisson process}.

A general \textit{pure jump-type homogeneous Markov process} on a measurable space $(E,\Sigma)$ can be defined
in the following way. Let $\lambda\colon E\to [0,\infty)$ be a given measurable function and let $P(x,B)$ be a
given \emph{transition probability} function on $E$,  i.e., $P(x,\cdot)$ is a probability measure for each $x\in
E$ and the function $x\mapsto P(x, B)$ is measurable for each $B\in\Sigma$. Let $t_0=0$ and let $X(0)=X_0$ be an
$E$-valued random variable.
 For each $n\ge 1$ we can choose the $n$th
\textit{jump time} $t_n$ as a positive random variable
satisfying
\[
\Prob(t_n-t_{n-1}\le t|X_{n-1}=x)=1-e^{-\lambda(x)t}, \quad t\ge 0,
\]
and we define
\[
X(t)=\left\{
  \begin{array}{ll}
X_{n-1} & \text{ for \,} t_{n-1}\le t<t_{n},\\
X_n &\text{ for \,} t=t_n,
\end{array}\right.
\]
where the $n$th \emph{post-jump position} $X_n$ is an $E$-valued random variable such that
\[
\Prob(X_n\in B|X_{n-1}=x)=P(x,B).
\]

Another type of simple PDMPs  is a \textit{velocity jump process}. An individual is moving in the space $\mathbb
R^d$ with a constant velocity  and at jump times  $(t_n)$ it chooses a new velocity. We assume that jump times
are the same as for a Poisson process $N(t)$ with intensity $\lambda$. It means that $F(t)=1-e^{-\lambda t}$ is
the probability distribution function of $t_n-t_{n-1}$. Let $x(t)$  be the position and $v(t)$ be the velocity
of an individual at time $t$.
 We assume that for every $x,v\in\mathbb R^d$,
there is
a probability Borel measure $P(x,v,B)$ on $\mathbb R^d$ which describes
the change of the velocity after a jump, i.e.,
\[
\Prob(v(t_n)\in B\,|x(t_n^{-})=x,\, v(t_n^{-})=v)=P(x,v,B)
\]
for every Borel subset $B$ of $\mathbb R^d$, where $x(t_n^{-})$ and $v(t_n^{-})$ are the left-hand side limits
of $x(t)$ and $v(t)$ at the point $t_n$. Between jumps
 the pair $(x(t),v(t))$ satisfies
the following system of ordinary differential equations
\begin{equation}
\left\{
\begin{aligned}
x'(t)&=v(t),\\
 v'(t)&=0.
\end{aligned}
\right.
\end{equation}
Then $X(t)=(x(t),v(t))$, $t\ge 0$, is a PDMP corresponding to this movement.

There are a number of interesting examples of velocity jump processes
 with applications to aggregation and chemotaxis
phenomena (see e.g. \cite{HH}). The simplest one is the symmetric movement on the real line $\mathbb R$. In this
case we assume that an individual is  moving with
 constant speed, say one, and at a jump time
it changes the direction of movement to the opposite one. A PDMP corresponding to the symmetric movement has
values in the space $\mathbb R\times\{-1,1\}$ and  $P(x,v,\{-v\})=1$ for $v=-1,1$. This process was first
studied by Goldstein \cite{Go} and  Kac \cite{Kac} in connection with the telegraph equation. It was called the
Goldstein-Kac \emph{telegraph process} afterwards and studied thoroughly in \cite{Ki}.

 More advanced  examples of velocity jump processes and their comparison  with  dispersal of cells,
insects and mammals are given in \cite{ODA,Stroock}. One can also consider velocity jump processes defined in a
bounded domain $G$. Examples of such processes are  stochastic billiards \cite{Evans} which do not change
velocity in the interior of $G$ but when an individual or a point hits the boundary, a new direction is chosen
randomly from directions that point back into the interior of $G$,  and the motion continues. PDMPs with jumps
at the boundary appear as well in the theory of gene regulatory systems, for example in a model of the
production of subtilin
 by the bacterium \textit{Bachillus subtilis} \cite{HWS}.

\subsection{Two phase cell cycle model}
\label{ss:tpcc} Now we consider another type  of  PDMPs which is a \textit{flow with jumps}  described in the
following way. Let $E$ be a topological space and let a continuous function $\pi\colon\mathbb{R}_+\times E\to E$
be a semiflow on $E$, i.e.,
\begin{enumerate}
\item[(a)] $\pi_0x=x$ for  $x\in E$,
\item[(b)]  $\pi_{s+t}x=\pi_t(\pi_sx)$ for  $x\in E$,
$s,t\in\mathbb{R}_+$.
\end{enumerate}
The semiflow $\pi_t$ describes the movement of points  between jumps,
 i.e.,  if $x$ is the position of the point at time $t$ then
 $\pi_{\tau}x$ is its position at time $t+\tau$.
The point located at $x$ can jump with an  intensity $\lambda(x)$ to a point $y$. The location of $y$ is
described by a transition function $P(x,B)$, i.e.,
 $P(x,B)$ is the probability that $y\in B$. After the jump it continues movement
according to the same principle.

A simple example of  a flow with jumps is the following size-structured model of a cellular population (see e.g.
\cite{mackeytyran08}). The cell size (mass, volume) $x>0$ grows with  rate $g(x)$ and
 it splits with intensity $\varphi(x)$ into two daughter cells with size $x/2$,
 i.e., $P(x,B)=1$ if $x/2\in B$ and   $P(x,B)=0$ otherwise.
After division we consider the size of a daughter cell, etc., and we obtain a process $X(t)$, $t>0$, which
describes the size of consecutive descendants of a single cell. The process $X(t)$, $t>0$, is a PDMP.

Another example of a flow with jumps  appears in the Rubinow model of a cellular population \cite{Rubinow}. In
this model we assume that a newborn cell has size $x=m$, then it grows with  rate $g(x)$ and when it reaches
size $x=2m$ it splits into two daughter cells with sizes $x=m$. Similarly to the previous model we consider a
process $X(t)$, $t>0$, which describes the size of consecutive descendants of a single cell. Although the jump
times in this process are  not random,  $X(t)$, $t>0$, is also a PDMP.

A more advanced flow with jumps is a two phase cell cycle model which is a combination of the two
size-structured models described above. The cell cycle is a series of events that
 take place in a cell leading to its replication \cite{MH}.
There are several models of the cell cycle  but from a mathematical point of view we can simplify these models
and we assume that there are only two phases in the cell cycle: the resting phase $A$ with a random duration
$t_A$, when a cell is growing, and the proliferating  phase $B$ with a constant duration $t_B$. Here we describe
a continuous time version of the Tyrcha model  \cite{Tyrcha} and we show that it can be treated as a PDMP. The
crucial role in the model  is played by a parameter $x$ which describes the state of a cell in the cell cycle.
It is not clear what $x$ exactly should be. We simply interpret $x$ as a cell size. The cell size $x>0$ grows
with rate $g(x)$ and the cell enters the phase $B$ with intensity $\varphi(x)$. It is clear that the process $X(t)$,
$t>0$, which describes the size of consecutive descendants of a single cell is piecewise deterministic but it is
non-Markovian because its  future $X(t)$, $t\ge t_0$, depends not only on the random variable $X(t_0)$ but also
on the phase in which it is at the time $t_0$.

Now we extend the process $X(t)$, $t\ge 0$, to obtain a homogeneous PDMP. A new process $\widetilde X(t)$, $t\ge
0$, is defined on the state space $[0,\infty)\times [0,t_B]\times \{1,2\}$ in the following way. Let $\widetilde
X(t)=(X(t),y,i)$, where $i=1$ if at time $t$ a cell is in the phase $A$ and $i=2$ if it is in the phase $B$.  We
let $y=0$ if the cell is in the phase $A$ and otherwise let $y$ be the time which elapsed since  the cell
entered the phase $B$. Let $s_n$ be a time when a cell from the $n$th generation enters the phase $B$. Since
the duration of the phase $B$ is constant and is equal to $t_B$, a cell from the $n$th generation splits at
time $t_n=s_n+t_B$. Between these jump points  the coordinates of the process $\widetilde X(t)$ satisfy the
following system of ordinary differential equations
\begin{equation}
\left\{
\begin{aligned}
\widetilde X'_1(t)&=g(\widetilde X_1(t)),\\
\widetilde X'_2(t)&=
\left\{
\begin{array}{ll}
0, \text{ \,\,if  $\widetilde X_3(t)=1$},\\
1, \text{ \,\,if  $\widetilde X_3(t)=2$},\\
\end{array}\right.
\\
\widetilde X'_3(t)&=0.
\end{aligned}
\right.
\end{equation}
The post-jump positions are given by
\[
\widetilde X_1(s_n)=\widetilde X_1(s_n^{-}),
\quad
\widetilde X_2(s_n)=\widetilde X_2(s_n^{-})=0,
\quad
\widetilde X_3(s_n)=2,
\]
and
\[
\widetilde X_1(t_n)=\tfrac12 \widetilde X_1(t_n^{-}),
\quad
\widetilde X_2(t_n)=0,
\quad
\widetilde X_3(t_n)=1.
\]
Let $\pi_tx_0=x(t)$ be the solution of the equation $x'=g(x)$ with initial condition $x(0)=x_0$.
The distribution function of $s_n-t_{n-1}$ is
given by
\begin{equation}
\label{tyr1} F(t)=1-\exp\Big\{- \int_0^{t}\varphi (\pi_s x_0) \,ds\Big\},
\end{equation}
where $x_0=\widetilde X_1(t_{n-1})$, while that of $t_n-s_n$ by $F(t)=0$ for $t< t_B$ and $F(t)=1$ for $t\ge
t_B$.

The life-span $t_n-t_{n-1}$ of a cell with initial size $x_0$ has the distribution function
\begin{equation}
\label{tyr1-cc}
F(t)=\left\{
\begin{array}{ll}
0, \text{ \,\,if  $t< t_B$},\\
1-\exp\Big\{- \int_0^{t-t_B}\varphi (\pi_s x_0) \,ds\Big\},
\text{ \,\,if  $t\ge t_B$},\\
\end{array}
\right.
\end{equation}
and we have the following relation between the distributions of the random variables $\widetilde{X}(t_{n})$ and
$\widetilde{X}(t_{n-1})$:
\begin{equation}
\label{tyr2cc} \widetilde{X}(t_{n})\stackrel{d}{=}\tfrac 12
\pi_{t_B}\Big(Q^{-1}\big(Q(\widetilde{X}(t_{n-1}))+\xi_n \big)\Big),
\end{equation}
where $\xi_n$ is a random variable independent of $\widetilde{X}(t_{n-1})$ with exponential distribution of mean
one and $ Q(x)=\int \limits_{0}^{x} \dfrac{\varphi (r)}{g(r)}\,dr $, $x>0$.

\subsection{Gene expression}
\label{ss:1-ge1} Another class of  PDMPs is the family of processes  with switching dynamics. Assume that we
have a finite number of semiflows $\pi^i_t$, $i\in I=\{1,\dots,k\}$
 on a topological space $E$.  The state of the system is a pair $(x,i)\in E\times I$.
 If the system is at state $(x,i)$ then $x$  can change according
 to the semiflow $\pi^i_t$ and after time $t$ reaches the state $(\pi^i_t(x),i)$
or it can switch to the state $(x,j)$ with  a bounded and continuous  intensity $q_{ji}(x)$. The pair
$(x(t),i(t))$ constitutes a Markov process $X(t)$  on $E\times I$.

Now we show how PDMPs can be applied to model gene expression. Gene expression is a complex process which
involves  three processes: gene activation/inactivation, mRNA transcription/decay,  and protein
translation/decay. We consider a simplified version of the  model of gene expression introduced by Lipniacki et
al. \cite{lipniacki2} and studied in \cite{BLPR}. We assume that the production of  proteins is regulated by a
single gene and we omit the intermediate process of mRNA transcription. A gene can be in an active or an
inactive state and it can be transformed into an active state or into an inactive state with intensities $q_0$
and $q_1$, respectively. The rates $q_0$ and $q_1$ depend on the number of protein molecules $X(t)$. If the gene
is active then proteins are produced with a constant speed $P$. Protein molecules undergo the process of
degradation with rate $\mu$ in both states of the gene. It means that the process $X(t)$, $t\ge 0$,  satisfies
the equation
\begin{equation}
\label{1:ga1}
X'(t)=P A(t) -\mu X(t),
\end{equation}
where $A(t)=1$ if the gene is active and $A(t)=0$ in the opposite case. Then the process $\widetilde
X(t)=(X(t),A(t))$, $t\ge 0$, is a PDMP. Since the right-hand side of equation (\ref{1:ga1}) is negative for
$X(t)>\frac P{\mu}$ we can restrict values of $X(t)$ to the interval $\big[0,\frac P{\mu}\big]$ and the process
$\widetilde X(t)$ is defined on the state space $\big[0,\frac P{\mu}\big]\times \{0,1\}$.

The process $\widetilde X(t)$ has jump points when the gene changes its activity. Formula (\ref{tyr1}) allows us
to find the distribution of the time between consecutive jumps. Observe that if $x_0$ is the number of protein
molecules at a jump time, then after time $t$ we have
\[
\pi^0_tx_0=x_0e^{-\mu t}, \quad \pi^1_tx_0=\frac P{\mu} +\Big(x_0-\frac P{\mu}\Big)e^{-\mu t},
 \]
 protein molecules in an inactive and an active state, respectively.
From (\ref{tyr1}) it follows that the probability distribution function of the length  of stay in an inactive state is
given by
\[
1-\exp\Big\{- \int_0^t q_0\Big (x_0e^{-\mu s}\Big ) \,ds\Big\}
\]
and in an active state by
\[
1-\exp\Big\{- \int_0^t q_0\Big(\frac P{\mu} +\Big(x_0-\frac P{\mu}\Big)e^{-\mu s}\Big ) \,ds\Big\}.
\]

\subsection{Neural activity}
\label{ss:neural}
A neuron is an electrically excitable cell that processes and transmits information through electrical signals.
 The neuron's membrane potential $V_m$ is the difference between the inside potential and the outside
potential. If a cell is in the resting state, then this potential, denoted by $V_{m,R}$,  is about ${}-70$ mV.
The \textit{depolarization} is defined as
\[
V=V_m - V_{m,R}.
\]
A cell is said to be \textit{excited} (or depolarized) if $V >0$ and \textit{inhibited} (or
hyperpolarized) if  $V < 0$.
The   Stein's model  \cite{Stein1965,Stein1967} describes how  the depolarization  $V(t)$ is changing in time.
The cell is initially at rest so that $V(0) = 0$.
Nerve cells may be excited or inhibited  through
neuron's synapses ---  junctions between nerve cells
(or between muscle and nerve cell) such that electrical activity in one cell may
influence the electrical potential in the other.  Synapses may be
excitatory or inhibitory.
We assume that there are two nonnegative constants $a_E$ and $a_I$
such that if at time $t$ an excitation  occurs then
$V(t^+)=V(t^-)+a_E$ and if
an inhibition occurs
then
$V(t^+)=V(t^-)-a_I$.
 The jumps (excitations and inhibitions)  may
occur at random times according to two independent Poisson processes $N^E(t)$, $N^I(t)$, $t\ge 0$, with positive
intensities  $\lambda_E$ and $\lambda_I$,  respectively. Between jumps  the depolarization $V(t)$ decays
according to the equation $V'(t)=-\alpha V(t)$. When a sufficient (threshold) level $\theta>0$ of excitation is
reached, the neuron emits an action potential (fires). This will  be followed by an absolute refractory period
of duration $t_R$, during which $V\equiv 0$ and then  the process starts again.

We now describe the neural activity as a PDMP. Since the refractory  period has a constant duration we can use a
model similar to that of  Section~\ref{ss:tpcc} with  two phases $A$ and $B$, where $A$ is the subthreshold
phase and $B$ is the refractory phase of duration $t_R$. We consider  two types of jump points: when  the neuron
is excited or inhibited and the ends of refractory periods. Thus, we can have one or more  jumps inside the
phase $A$.

 Let $\widetilde X(t)=(V(t),0,1)$ if the neuron is in the phase $A$ and
$\widetilde X(t)=(V(t),y,2)$ if the neuron is in the phase $B$, where $y$ is the time since the moment of
firing. The process $\widetilde X(t)=(\widetilde X_1(t),\widetilde X_2(t),\widetilde X_3(t))$, $t\ge 0$, is
defined on the state space $(-\infty,\theta)\times [0,t_R]\times \{1,2\}$ and between jumps it satisfies the
following system of equations
\begin{equation}
\left\{
\begin{aligned}
\widetilde X'_1(t)&=- \alpha \widetilde X_1(t),\\
\widetilde X'_2(t)&=
\left\{
\begin{array}{ll}
0, \text{ \,\,if  $\widetilde X_3(t)=1$},\\
1, \text{ \,\,if  $\widetilde X_3(t)=2$},\\
\end{array}\right.
\\
\widetilde X'_3(t)&=0.
\end{aligned}
\right.
\end{equation}
Let $t_0,t_1,t_2,\dots$ be the subsequent jump times. We denote by $\mathcal F$ the subset of jump times
consisting of  firing points. If the neuron is in the phase $A$, i.e.,  $\widetilde X_3(t)=1$, the
depolarization can jump with intensity $\lambda=\lambda_E+\lambda_I$. It means that $F(t)=1-e^{-\lambda t}$ is
the distribution function of $t_n-t_{n-1}$ if $t_{n-1}\notin \mathcal F$. If $t_{n-1}\in \mathcal F$ then the
distribution of $t_n-t_{n-1}$ is  $F(t)=0$ for $t< t_R$ and $F(t)=1$ for $t\ge  t_R$. The transition at a jump
point depends on the state of the neuron (its phase and the value of its depolarization). If $\widetilde
X(t_n^-)=(0,t_R,2)$ then
 $\widetilde X(t_n)=(0,0,1)$ with probability one;
 if $\widetilde X(t_n^-)=(x,0,1)$ and  $x<\theta-a_E$ then
 $\widetilde X(t_n)=
(x+a_E,0,1)$ with probability $\lambda_E/\lambda$
and
 $\widetilde X(t_n)=
(x-a_I,0,1)$ with probability $\lambda_I/\lambda$;
while
 if $\widetilde X(t_n^-)=(x,0,1)$ and  $x\ge \theta-a_E$ then
 $\widetilde X(t_n)=
(0,0,2)$ with probability $\lambda_E/\lambda$ and
 $\widetilde X(t_n)=
(x-a_I,0,1)$ with probability $\lambda_I/\lambda$.

\subsection{Size-structured population  model}
\label{ss:ppr-psm} In this section we return  to size-structured models but instead of a single cell line  we
consider the size distribution of  all cells in the population. This model can serve as a prototype of
individual based models like age and phenotype structured models as well as models of coagulation-fragmentation
processes.

The size $x(t)$ of a cell grows according
to the equation
\[
x'(t)=g(x(t)).
\]
A  single cell with size $x$ replicates with rate $b(x)$  and dies with rate $d(x)$. A daughter cell has a
half size of the mother cell. Let us assume that at time $t$ we have $k$ cells and denote by
$x_1(t),x_1(t),\dots,x_k(t)$ their sizes. We can assume that  a state of the population  at time $t$ is the set
\[
\{x_1(t),\dots,x_k(t)\}
\]
and that the evolution of the population is  a stochastic process $$X(t)=\{x_1(t),\dots,x_k(t)\}.$$ Since the
values of this process are sets of points the process $X(t)$ is called a \textit{point process}. Thought such
approach is a natural one it has one important disadvantage.  We are unable to describe properly a situation
when two cells have the same size. One solution of this problem is to consider $X(t)$  as a process whose values
are  multisets. We recall that a \textit{multiset}\index{multiset} (or a \textit{bag}) is a generalization of
the notion of a set in which members are allowed to appear more than once. Another artful solution of this
problem is to consider $X(t)$ as a process with values in the space of measures given by
\[
X(t)=\delta_{x_1(t)}+\dots+\delta_{x_k(t)},
\]
where $\delta_a$ denotes the \textit{Dirac measure}\index{Dirac measure} at point $a$, i.e., $\delta_a$ is the
probability measure concentrated at the point $a$. This approach has some disadvantages also, for example it is
rather difficult to consider differential equations on measures. Yet another solution of this problem is to
consider a state of the system as $k$-tuples   $(x_1(t),\dots,x_k(t))$. Since some cells can die or split into
two cells, the length of the tuple changes in time. To omit this difficulty we introduce an extra "death state"
$*$ and we describe the state of the population at time $t$ as an infinite sequence of elements from the space
$\mathbb R^+_*=\mathbb [0,\infty)\cup\{*\}$ which has numbers  $x_1(t),\dots,x_k(t)$ on some $k$ positions and
it has $*$ on the remaining positions. In order to have uniqueness of states we introduce  an equivalence
relation $\sim$ in the space  $E$ of all $\mathbb R^+_*$- valued sequences $x$ such that $x_i=*$ for all but
finitely many $i$. Two sequences $x\in E$ and $y\in E$ are equivalent with respect to $\sim$ if $ y$  can be
obtained as a permutation of $x$, i.e., $ x\sim y$ if and only if there is a bijective function $\sigma:\mathbb
N\to \mathbb N$ such  that $y=(x_{\sigma(1)},x_{\sigma(2)},\dots)$. The state space $\widetilde E$ in our model
is the space of all equivalence classes with respect to $\sim$, i.e., $\widetilde E=E/\sim$.

Now we can describe the evolution of the population as a stochastic process  $X(t)=[(x_1(t),x_2(t),\dots)] $
with values in the space $\widetilde E$ where $[ ]$ denotes an equivalence class. The process $X(t)$ has jump
points when one of the cells dies or replicates. We define  $g(*)=b(*)=d(*)=0$ and  admit the convention that
$x(t)=*$ is the solution of the equation $x'(t)=0$ with initial condition $x(0)=*$. Between jumps the process
$X(t)$ satisfies the equation
\begin{equation}
\label{1:grow-ssp}
[(x_1'(t)-g(x_1(t)),x_2'(t)-g(x_2(t)),\dots)]=[(0,0,\dots)].
\end{equation}

For $t\ge 0$ and $x^0\in \mathbb R^+_*$ we denote by
$\pi(t,x^0)$
the solution $x(t)$ of the equation $x'(t)=g(x(t))$ with initial condition $x(0)=x_0$.
Let $\mathbf x^0=[(x_1^0,x_2^0,\dots)]\in \widetilde E$ and define
\[
\tilde \pi_t\mathbf x^0=[(\pi_t x_1^0,\pi_tx_2^0,\dots)].
\]
The jump rate function  $\varphi(\mathbf x)$
 at state $\mathbf x=[(x_1,x_2,\dots)]$ is the sum of rates of deaths and divisions of all cells:
\begin{equation}
\label{1:grow-inten}
\varphi(\mathbf x)=\sum_{i=1}^{\infty} (b(x_i)+d(x_i)).
\end{equation}
If $\mathbf x^0 \in \widetilde E$ is the initial state of the population at a jump time $t_n$, then
the probability distribution function
of $t_{n+1}-t_n$
is given by
\begin{equation}
\label{1:time} 1-\exp\Big\{- \int_0^t \varphi(\tilde \pi_s\mathbf x^0) \,ds\Big\}.
\end{equation}
At time $t_n$ one of the cells dies or replicates. If a cell dies we change the sequence by removing the cell's
size from the sequence  and we have
\begin{equation*}
\Prob\big(X(t_n)=[(x_1(t_n^-),\dots,x_{i-1}(t_n^-),x_{i+1}(t_n^-),\dots)]\big) =\frac
{d_i(x_i(t_n^-))}{\varphi(X(t_n^-))}
 \end{equation*}
for $i\in \mathbb N$. If a cell replicates  we remove its size from the sequence and add two new elements in the
sequence with sizes of the daughter cells and we have
\begin{multline*}
\Prob\big(X(t_n)=[(x_1(t_n^-),\dots,x_{i-1}(t_n^-),\tfrac12 x_i(t_n^-),\tfrac12
x_i(t_n^-),x_{i+1}(t_n^-),\dots)]\big)\\ =\frac {b_i(x_i(t_n^-))}{\varphi(X(t_n^-))}
 \end{multline*}
for $i\in \mathbb N$. In this way we have checked that the point process $X(t)$, $t\ge 0$, is a homogeneous PDMP
with values in $\widetilde E$.

We can identify the space $\widetilde E$ with the space $\mathcal N$ of finite counting measures on $\mathbb{R}_+$
by a map $\eta\colon \widetilde E\to \mathcal N$ given by
\begin{equation}
\label{1:izom}
\eta(\mathbf x)=  \sum_{\{i:\,\,x_i\ne *\}} \delta_{x_i}
\end{equation}
where $\mathbf x=[(x_1,x_2,\dots)]$. It means that the process $\eta(X(t))$, $t\ge 0$,  is a homogeneous PDMP
with values in $\mathcal N$.

\begin{remark}
\label{1:r1} In order to describe the jump transformation at times $t_n$ we need, formally,  to introduce a
$\sigma$-algebra $\Sigma$ of subset of $\widetilde E$ to define a transition function $P: \tilde E\times \Sigma
\to [0,1]$. Usually,  $\Sigma$ is a $\sigma$-algebra of Borel subsets of $\widetilde E$, thus we need to
introduce a topology on the space $\widetilde E$. Since  the space $\mathcal N$ is equipped with the topology of
weak convergence of measures, we can define open sets in $\widetilde E$ as preimages  through the function
$\eta$ of open sets in $\mathcal N$. Another way to introduce a topology is to construct directly a metric on
the space $\widetilde E$. Generally, a point process describes the evolution of configurations of points in a
state space which is a metric space $(S,\rho)$. First, we extend the state space $S$ by adding "the death
element" $*$. We need to define a metric on $S\cup \{*\}$. The best situation is if $S$ is a proper subset of a
larger metric space $S'$. Then we simply choose $*$ as an element from $S'$ which does not belong to the closure of $S$ and we
keep the same metric. In the other case, first we choose $x_0\in S$ and define $\rho(*,x)=1+\rho(x_0,x)$ for
$x\in S$. Next, we define a metric $d$ on the space $E$ by
\[
d(x,y)=\max_{i\in\mathbb N} \rho(x_i,y_i)
\]
and, finally,
we define
a metric $\widetilde d$ on the space $\widetilde E$ by
\[
\widetilde d([x],[y])= \min\{d(a,b)\colon \,\, a\in [x], \,\, b\in [y]\}.
\]
We next show that the topology in $\widetilde E$ induced from $\mathcal N$ is equivalent to the topology defined
by $\widetilde d$. Indeed, a sequence $(\mu_n)$ of finite counting measures converges weakly to a finite
counting measure $\mu$ iff the measures $\mu$ and $\mu_n$, $n\ge 1$ can be represented in the form
\[
\mu= \sum_{i=1}^{k}\delta_{x_{i}}, \quad \mu_n= \sum_{i=1}^{k_n}\delta_{x_{i,n}},
\]
where $k_n=k$ for sufficiently large $n$ and $\lim_{n\to\infty} \rho(x_{i,n},x_i)=0$ for $i=1,\dots,k$. Thus the
convergence of counting measures implies that the sequence $x^n=(x_{1,n},\dots,x_{k_n,n})$ converges to
$x=(x_1,\dots,x_k)$ in the metric $d$, and the sequence $[x^n]$  converges  to $[x]$ in $\widetilde d$. The
proof of the opposite implication goes in the same way.
\end{remark}

\section{Stochastic semigroups}
\label{s:stoch}
Most of PDMPs define stochastic semigroups  which describe the evolution of
 densities of the distribution of these processes.
In this section we recall the definition of a stochastic semigroup and provide a couple of examples of such
semigroups.

Let the triple $(E,\Sigma,m)$ be a $\sigma$-finite measure space.
Denote by $D$ the subset of the space
$L^1=L^1(E,\Sigma,m)$ which contains all
densities
\[
D=\{f\in \, L^1: \,\, f\ge 0,\,\, \|f\|=1\}.
\]
A linear operator $P\colon  L^1\to L^1$ is called a
\textit{stochastic (or Markov) operator\,}
 if $P(D)\subset D$.
Let $\{P(t)\}_{t\ge0}$ be a $C_0$-\textit{semigroup}, i.e.,
it satisfies the following conditions:
\begin{itemize}
\item[(a)] \  $P(0)=I$,  i.e., $P(0)f =f$,
\item[(b)] \ $P(t+s)=P(t) P(s)\quad \textrm{for}\quad
s,\,t\ge0$,
\item[(c)] \ for each $f\in L^1$ the function
$t\mapsto P(t)f$ is continuous.
\end{itemize}
Then the $C_0$-semigroup  $\{P(t)\}_{t\ge0}$  is called \textit{stochastic} iff each operator $P(t)$ is
stochastic. The \textit{infinitesimal generator} of $\{P(t)\}_{t\ge0}$ is  the operator $A$ with domain
$\mathcal{D}(A)\subseteq L^1$ defined as
\[
 A
f=\lim_{t\downarrow 0}\frac{1}{t}(P(t)f-f),\quad\mathcal{D}(A)=\{f\in L^1: \lim_{t\downarrow
0}\frac{1}{t}(P(t)f-f) \text{ exists}\}.
\]

Our first example of a stochastic semigroup is the following. Let $g\colon \mathbb R^d\to \mathbb R^d$ be a
$C^1$ function and consider the differential equation
  \begin{equation}
\label{M1}
x' (t) = g(x (t) ).
\end{equation}
Assume that $E$ is a measurable subset of $\mathbb R^d$
with a positive Lebesgue measure
such that
for each point $x_0\in E$ the solution $x(t)$ of  (\ref{M1}) with
$x(0)=x_0$ exists  and $x(t)\in E$ for all $t\ge 0$.
We denote this solution by $\pi_tx_0$. Let
$\Sigma$ be the $\sigma$-algebra of the Borel subsets of $E$ and
$m$ be the Lebesgue on $E$. Let $f\colon E\to [0,\infty)$ be a density
 and let $X_0$ be a random vector with values in $E$ with density $f$, i.e.,
 $\Prob(X_0\in B)=\int_B f(x)\,dx$ for  each Borel subset $B$ of $E$.
 Let $X(t)=\pi_t X_0$. Then the density of the distribution of
the  random vector $X(t)$
is given by
 \[
P(t)f(x)=
\begin{cases}
f(\pi_{-t}x) \det\Big[ \dfrac{d}{d x} \pi_{-t}x\Big], \quad
\textrm{if  $x\in \pi_t(E)$,}\\
0,\quad \textrm{if  $x\notin \pi_t(E)$,}
\end{cases}
\]
where $\pi_{-t}$ denotes the inverse of the one-to-one and onto mapping $\pi_{t}\colon E\to \pi_t(E)$. The
operators $P(t)$, extended linearly from $D$ to $L^1$, form a stochastic semigroup. If $f$ is a $C^1$ function
then the function $u(t,x)=P(t)f(x)$ satisfies the following partial differential equation
 \begin{equation}
\label{M2}
\frac{ \partial u (t,x) } {\partial t} =
- {\rm div} (g(x) u(t,x) ).
\end{equation}
If $A$ is an infinitesimal generator of the
semigroup $\{P(t)\}_{t\ge0}$ then
\begin{equation}
\label{M3} Af(x)= -{\rm div} (g (x) f(x)) = -\sum_{i=1}^d \frac{\partial}{\partial x_i} (g_i(x)f(x) ).
\end{equation}

 Now we consider
the processes $X(t)=(x(t),i(t))$   with switching dynamics described in Section~\ref{ss:1-ge1}. We assume that
each flow $\pi^i_t$,  $i\in I=\{1,\dots,k\}$, is defined as the solution of a system of differential equations
$x'=g_i(x)$ on a measurable subset $E$ of $R^d$. Let   $\{S^i(t)\}_{t\ge0}$ be the stochastic semigroup related
to $\pi^i_t$ and let the operator $A_i$  be its generator. If $f=(f_1,\dots,f_k)$ is a vertical vector
consisting of functions $f_i$ such that $f_i\in \mathcal D(A_i)$, we set
 $Af=(A_1f_1,\dots, A_kf_k)$ which is also a vertical vector.
We define $q_{jj}(x)=-\sum_{i\ne j}q_{ij}(x)$ and denote by $Q(x)$ the matrix $[q_{ij}(x)]$. Then the process
$X(t)$ induces a stochastic semigroup on the space $L^1(E\times I,\mathcal B(E\times I),m)$ with the
infinitesimal generator $Q+A$. Here
 $\mathcal B(E\times I)$ is the
$\sigma$-algebra of Borel
subsets of  $E\times I$ and
 $m$ is the product
measure on $\mathcal B(E\times I)$ given by $m(B\times\{i\})=\mu (B)$.

Finally, we provide stochastic semigroups for the flows with jumps $X(t)$ from Section~\ref{ss:tpcc}. Let
$\pi_tx$ be the semiflow describing solutions of equation~\eqref{M1} and let $\lambda(x)$ be the intensity of
jumping from the point $x$ to a point $y\in B$ chosen according to the transition probability $P(x,B)$. Suppose
that there is a stochastic operator $P$ on $L^1(E,\Sigma,m)$ induced by $P(x,\cdot)$, i.e.,
\begin{equation}\label{d:top}
\int_{E} P(x,B)f(x)m(dx)=\int_B Pf(x)m(dx)\quad \text{for all } B\in \Sigma, f\in D.
\end{equation}
If $\lambda$ is bounded then the process $X(t)$ induces a stochastic semigroup on the space $L^1(E,\Sigma,m)$
with infinitesimal generator of the form $A_0f-\lambda f+P(\lambda f)$, where $A_0$ is as $A$ in \eqref{M3}. If
$\lambda$ is unbounded then one may need to impose additional constraints on $A_0$, $\lambda$, and/or $P$ to
obtain a stochastic semigroup for $X(t)$, see \cite{tyran09,tyran09b} for necessary and sufficient conditions.
For the particular example of the model of the cell cycle on $E=(0,\infty)$ with  one phase we have
$Pf(x)=2f(2x)$ for $x>0$. Suppose that $g\colon (0,\infty)\to (0,\infty)$ is continuous and $\varphi/g$ is
locally integrable with
\begin{equation}
\int_{\bar{x}}^{\infty} \frac{1}{g(r)}dr=\int_{\bar{x}}^{\infty} \frac{\varphi(r)}{g(r)}dr=\infty
\end{equation}
for some $\bar{x}>0$. Then the process $X(t)$ induces a stochastic semigroup $\{P(t)\}_{t\ge 0}$ on
$L^1=L^1((0,\infty),\Sigma,m)$, where $\Sigma$ is the Borel $\sigma$-algebra of subsets of $(0,\infty)$ and $m$
is the Lebesque measure, with infinitesimal generator of the form \cite{mackeytyran08}
\[
Af(x)=-\frac{d }{d x}(g(x)f(x))-\varphi(x)f(x)+2\varphi(2x)f(2x)
\]
defined for $f\in\mathcal{D}(A)=\mathcal{D}_0\cap L_\varphi^1$, where $L_\varphi^1=\{f\in L^1:
 \varphi f\in L^1\}$ and
\[
\mathcal{D}_0=\{f\in L^1: gf\text{ is absolutely continuous}, \;(gf)'\in L^1\},
\]
together with the boundary condition $\lim_{x\to 0}g(x)f(x)=0$. 

 For the two phase model we can restrict the state space to the set $(0,\infty)\times\{0\}\times\{1\}\cup
(0,\infty)\times[0,t_B]\times\{2\}$. We consider the corresponding stochastic semigroup $\{P(t)\}_{t\ge 0}$ on
the product space $L^1((0,\infty))\times L^1((0,\infty)\times [0,t_B])$. Let $f=(f_1,f_2)$
be the  density of the process at time $t$, where
$f_1(t,x)$ and $f_2(t,x,y)$ denote the
partial densities related to the phases $A$ and $B$, respectively.
If $f_1$, $f_2$ are smooth functions  then they satisfy the following
equations
\[
\begin{split}
\frac{\partial f_1(t,x)}{\partial t}&=-\frac{\partial }{\partial
x}(g(x)f_1(t,x))-\varphi(x)f_1(t,x)+2f_2(t,2x,t_B),\\
\frac{\partial f_2(t,x,y)}{\partial t}&=-\frac{\partial }{\partial x}(g(x)f_2(t,x,y))-\frac{\partial }{\partial
y}(f_2(t,x,y)),
\end{split}
\]
with the boundary conditions
\[
f_2(t,x,0)=\varphi(x)f_1(t,x),\quad x>0,\, t\ge 0,
\]
\[
\lim_{x\to 0} g(x)f_1(t,x)=\lim_{x\to 0} g(x)f_2(t,x,y)=0,\quad y\in [0,t_B],\, t\ge 0.
\]

\section{Long time behavior}
\label{s:ltb} In this section we study asymptotic properties of stochastic semigroups induced by PDMPs. We will
consider two properties: asymptotic stability and sweeping. A stochastic semigroup $\{P(t)\}_{t\ge 0}$ on
$L^1(E,\Sigma,m)$ is called {\it asymptotically stable\,}
 if there is a density
$f_*$ such that
\begin{equation}
\label{as-def}
\lim _{t\to\infty}\|P(t)f-f_*\|=0 \quad \text{for}\quad
f\in D.
\end{equation}
A density $f_*$ which satisfies (\ref{as-def}) is \textit{invariant\,}, i.e.,  $P(t)f_*=f_*$ for each $t>0$.
A~stochastic semigroup $\{P(t)\}_{t\ge 0}$ is called \textit{sweeping} with respect to a set $B\in\Sigma$ if for
every
 $f\in D$
\begin{equation*}
\lim_{t\to\infty}\int_B P(t)f(x)\,m(dx)=0.
\end{equation*}

Let us now recall two general results concerning asymptotic properties of partially integral semigroups. A
stochastic semigroup $\{P(t)\}_{t\ge 0}$ on $L^1(E,\Sigma,m)$ is called {\it partially integral} if there exists
a measurable function $k\colon (0,\infty)\times E\times E\to[0,\infty)$, called a {\it kernel}, such that
\begin{equation*}
P(t)f(x)\ge\int_E k(t,x,y)f(y)\,m (dy)
\end{equation*}
for every density $f$ and
\begin{equation*}
\int_E\int_E  k(t,x,y)\,m(dy)\,m(dx)>0
\end{equation*}
for some $t>0$.
 \begin{theorem}[\cite{PR-jmaa2}]
\label{asym-th2}
Let $\{P(t)\}_{t\ge 0}$ be a partially integral stochastic
semigroup. Assume that the  semigroup $\{P(t)\}_{t\ge 0}$ has
a unique invariant density $f_*$. If $f_*>0$ a.e., then the semigroup
$\{P(t)\}_{t\ge 0}$ is asymptotically stable.
\end{theorem}

To prove asymptotic stability,  it is  sometimes difficult to check directly that the  semigroup $\{P(t)\}_{t\ge
0}$ has a unique invariant density $f_*$. Therefore, the following theorem can be useful in checking whether a
semigroup is asymptotically stable or sweeping.

\begin{theorem}[\cite{R-b95}]
\label{KTR}
Let $E$ be a metric space and $\Sigma=\mathcal B(E)$ be the
$\sigma$--algebra of Borel subsets of $E$. We assume  that a partially integral  stochastic
semigroup $\{P(t)\}_{t\ge 0}$ with the kernel $k$ has the following properties:
\item[\,\,(a)]
 for every $f\in D$ we have
$\int_0^{\infty}P(t)f\,dt>0$ a.e.,
\item[\,\,(b)]
for every $y_0\in E$ there exist $\varepsilon >0$, $t>0$,
 and a measurable  function
$\eta \ge 0$ such that $\int \eta\, dm >0$ and
\begin{equation*}
k(t,x,y)\ge \eta (x)
\end{equation*}
for $x\in E$ and $y\in B(y_0,\varepsilon)$, where
$B(y_0,\varepsilon)$ is the open ball with center $y_0$ and radius
$\varepsilon$.
Then the semigroup $\{P(t)\}_{t\ge 0}$
is asymptotically stable if it
 has an invariant density and it is sweeping
with respect to compact sets if it
has no invariant density.
In particular, if $E$ is a compact set  then
the semigroup $\{P(t)\}_{t\ge 0}$
is asymptotically stable.
\end{theorem}

We are now ready to apply Theorems~\ref{asym-th2} and \ref{KTR} to  stochastic semigroups induced by PDMPs with
switching dynamics. In many applications a PDMP with switching dynamics is induced by flows $\pi_t^i$, $i\in
I=\{1,\dots,k\}$, acting on an open subset $G$ of $\mathbb R^d$, and we start with a stochastic semigroup
$\{P(t)\}_{t\ge 0}$ defined on the space $L^1(G\times I,\mathcal B(G\times I),dx\times di)$, but this semigroup
has  a stochastic attractor having some additional properties. By a \textit{stochastic attractor} we understand
here a measurable subset $S$ of $G$ such that for every density $f\in L^1(G\times I)$ we have
\begin{equation}
\label{st-at}
\lim_{t\to\infty}\int\limits_{S\times I}
P(t)f(x,i)\,dx\,di=1.
\end{equation}
For example, if there exists a measurable subset $S$ of $G$ such that $x(t)(\omega)\in S$ for $t>t(\omega)$ for
almost every $\omega$ then $S$ is a stochastic attractor. If a stochastic semigroup has a stochastic attractor
$S$ then it is enough to study the restriction of the semigroup $\{P(t)\}_{t\ge 0}$ to the space $L^1(E,\mathcal
B(E),m)$, where $E=S\times I$ and $dm=dx\times di$.

Let us now explain  how to check conditions (a) and (b) of Theorem~\ref{KTR}.  We can obtain condition (a) if we
check that  $m$-almost every two states $(x,i)\in E$, $(y,j)\in E$ can be joined by
 a path  of the process $(x(t),i(t))$. To be precise
there exist $n\in\mathbb N$, $\mathbf i=(i_1,\dots,i_n)\in I^n$,  and $\mathbf t=(t_1,\dots,t_n)\in
(0,\infty)^n$ such that $i_1=i$, $i_n=j$, and
\[
y=\mathbf \pi^{\mathbf i}_{\mathbf t}(x)=
\pi
_{t_n}^{i_{n}}
\circ\cdots\circ \pi
_{t_2}^{i_2}
\circ
\pi
_{t_1}^{i_1}
(x).
 \]

Condition (b) can be checked by using Lie brackets. We now recall the definition of Lie brackets. Let $a(x)$ and
$b(x)$ be two vector fields on $\mathbb R^d$. The \textit{Lie bracket} $[a,b]$ is a vector field given by
\[
[a,b]_j(x)=
\sum_{k=1}^d\left(
 a_k
\frac{\partial b_j}{\partial x_k}(x)
- b_k
\frac{\partial a_j}{\partial x_k}(x)
\right).
\]
Let a PDMP with switching dynamics be defined by  the systems  of differential equations $x'=g_i(x)$, $i\in I$,
with intensities $q_{ji}(x)$. We say that the {\it H\"{o}rmander's condition} holds at a point $x$ if $q_{ij}
(x)>0$ for  all $1\le i, j\le k$ and if  vectors
\[
g_2(x)-g_1(x), \dots, g_k(x)-g_1(x),\, [g_i,g_j](x)_{1\leq i,j\leq k},\, [g_i,[g_j,g_l]](x)_{1\leq i,j,l\leq
k},\dots
\]
span the space $\mathbb R^d$. Let $y_0\in S$ and assume that there exist $n\in\mathbb N$, $\mathbf i\in I^n$ and
$\mathbf t\in (0,\infty)^n$ such that the H\"{o}rmander's condition holds at the point $y$ given by
\[
y=\mathbf \pi^{\mathbf i}_{\mathbf t}(y_0),
\]
then $y_0$ satisfies condition (b). This fact is a simple consequence of \cite[Theorem 4]{BH}.

Finally, we give some examples of applications to biological models.

\begin{example}[Gene expression]
\label{ex:gene expression} The model of gene expression from Section~\ref{ss:1-ge1} is a special case of the
following PDMP with switching dynamics. We have two flows induced by one-dimensional differential equations
$x'=g_i(x)$, $i=1,2$, where $g_i$ are $C^{\infty}$-functions with the following property: there exist points
$x_1<x_2$ such that
\[
g_i(x)>0\textrm{ for $x<x_i$ and }
g_i(x)<0\textrm{ for $x>x_i$.}
\]
It is obvious that almost all trajectories enter  the set $S=[x_1,x_2]$. Observe that any two states $(x,i)$ and
$(y,j)$ with $x,y\in (x_1,x_2)$ and $i,j\in \{0,1\}$ can be joined by
 a path  of the process $(x(t),i(t))$. Hence, condition (a) of  Theorem~\ref{KTR} is fulfilled. Since
$g_2(x)-g_1(x)>0$ for $x\in S$, the H\"{o}rmander's condition holds at each point $x\in S$ and, therefore,
condition (b)  is fulfilled. Since the set $E=S\times \{1,2\}$ is compact, the semigroup induced by our PDMP is
asymptotically stable. More precisely, there exists a density $f_*\colon \mathbb R\times\{1,2\}\to [0,\infty)$
such that $f_*(x,i)=0$ for $x\notin [x_1,x_2]$ and $\lim_{t\to\infty} \|P(t)f-f_*\|=0$ for every density $f\in
L^1(\mathbb R\times\{1,2\})$.
\end{example}

\begin{example}[Population model with and without Allee effect]
\label{ex:Allee efect}
Consider a PDMP
with switching dynamics
induced by two differential equations
\begin{equation}
\label{Allee1}
x' (t) =\lambda  \bigg(1- \frac{x(t)}{K}- \frac{Ai}{1+Bx(t)} \bigg) x (t),
\end{equation}
where $i=0,1$ and $A,B,K$ are positive constants such that $KB>1$ and
\begin{equation}
\label{Allee2}
1<A<\frac{(BK+1)^2}{4KB}.
\end{equation}
The number $x(t)>0$ describes the size of a population. If $i=0$ then (\ref{Allee1}) reduces to a logistic
equation and $\lim_{t\to\infty} x(t)=K$. If $i=1$ then  (\ref{Allee1}) has three stationary states $x_0,x_1,x_2$
such that $x_0=0<x_1<x_2<K$ with the following properties. If $x(0)\in (0,x_1)$ then $\lim_{t\to\infty} x(t)=0$
(called Allee effect) and if $x(0)\in (x_1,\infty)$ then $\lim_{t\to\infty} x(t)=x_2$. Now we consider a PDMP
induced by these equations with positive and continuous intensities of switching. Almost all trajectories enter
the interval $S=[x_2,K]$, thus $S$ is a stochastic attractor. Since almost all states in $E=S\times\{0,1\}$ can
be joined by
 paths  of the process $(x(t),i(t))$  and $g_0(x)>g_1(x)$, the assumptions of   Theorem~\ref{KTR}  are fulfilled and the semigroup
induced by our process is asymptotically stable.
  \end{example}

\begin{example}[Population model with two different birth rates]
\label{ex:birth} Now we consider a population model with a constant death rate $\mu$ and birth rates
$b_i(x)=b_i-cx$, $i=0,1$, which can change in time. Thus, the size $x\ge 0$ of the population is described by  a
PDMP with switching dynamics defined by  two differential equations
\[
x'=g_i(x)
\]
with $g_i(x)= (b_i-cx)x-\mu x$ for  $i=0,1$. Denote by $q_i(x)$  the intensities of changing the state $i$ to
$1-i$. We assume that $b_0<\mu$ and $b_1>\mu$ and that the intensities $q_i(x)$ are continuous, positive, and
bounded functions. Observe that $g_i(0)=0$ for $i=0,1$, $g_0(x)<0$  for $x>0$ and that there exists a point
$a>0$ such that $g_1(x)>0$ for $x\in (0,a)$ and $g_1(x)<0$ for $x>a$.   The interval $S=(0,a]$ is a stochastic
attractor for this PDMP. Since almost all states in $E=S\times\{0,1\}$ can be joined by
 paths  of the process $(x(t),i(t))$ and $g_0(x)<g_1(x)$ for $x>0$, conditions (a) and (b) of  Theorem~\ref{KTR}
are fulfilled. Consequently, the semigroup  induced by our process is asymptotically stable or sweeping from
compact subsets of $E$.

In order to get asymptotic stability of  this semigroup, we need to check whether this
semigroup has an invariant density. Observe that if  $f(x,i)$  is
 an invariant density then the functions $f_i(x)=f(x,i)$ for $i=0,1$
 should be stationary solutions of the Fokker-Planck equation, i.e., $f_0$, $f_1$   satisfy
the following system of differential equations
\begin{equation}
\label{st-FP}
\left\{
\begin{aligned}
(g_0(x)f_0(x))'&=q_1(x)f_1(x)-q_0(x)f_0(x),\\
(g_1(x)f_1(x))'&=q_0(x)f_0(x)-q_1(x)f_1(x).
\end{aligned}
\right.
\end{equation}
Fix a point $x_0\in (0,a)$ and let
\[
r(x)=\frac{q_0(x)}{g_0(x)}+\frac{q_1(x)}{g_1(x)}\text{ and }
R(x)=\int_{x_0}^x r(s)\,ds.
\]
Then the functions
\[
\bar f_0(x)= -e^{-R(x)}/g_0(x) \text{ and } \bar f_1(x)= e^{-R(x)}/g_1(x)
\]
are positive in the interval $(0,a)$ and  they satisfy the system (\ref{st-FP}). If
\begin{equation}
\label{nier}
\alpha =\int_0^a (\bar f_0(x)+\bar f_1(x))\,dx<\infty,
\end{equation}
then  the semigroup  $\{P(t)\}_{t\ge 0}$ has an invariant density $f_*(x,i)$  given by $f_*(x,i)=\alpha^{-1}\bar
f_i(x)$, $i=0,1$, and, consequently, this semigroup is asymptotically stable.

If inequality (\ref{nier}) does
not hold, then the semigroup  $\{P(t)\}_{t\ge 0}$ has no invariant density. Indeed, if it has an invariant
density, say $f_\diamond (x,i)$, then the semigroup $\{P(t)\}_{t\ge 0}$ should be asymptotically stable and, at
the same time, if we extend the semigroup $\{P(t)\}_{t\ge 0}$
  to nonnegative measurable functions then $\bar f(x,i)$ is a non-integrable
stationary point of this semigroup.
Let us define $f^n(x,i)=f_i^n(x)=\bar f_i(x)\wedge n$  and
$c_n=\int_0^a (f^n_0(x)+f^n_1(x))\,dx$.
Then $\lim_{n\to\infty} c_n=\infty$ and
\[
\liminf_{t\to \infty} P(t)\bar f\ge  \lim_{t\to \infty} P(t)f^n= c_nf_\diamond,
\]
for all $n\in\mathbb N$, which contradicts the fact that $P(t)\bar f=\bar f$. Thus,  if  inequality (\ref{nier})
does not hold, then the semigroup  $\{P(t)\}_{t\ge 0}$ has no invariant density and according to
Theorem~\ref{KTR} this semigroup is sweeping from compact subsets of $E$. Since $[\varepsilon,a]\times \{0,1\}$
is a compact subset of $E$ for each $\varepsilon\in (0,a)$, sweeping means here that
\begin{equation}
\label{ex:sweep}
\lim_{t\to\infty} \int\limits_{[0,\varepsilon]\times\{0,1\}} P(t)f(x,i)\,m(dx,di)=1.
\end{equation}
Let  $p_i=q_{1-i}(0)/(q_0(0)+q_1(0))$. The numbers $p_i$ can be interpreted as the mean time of staying in the
state $i$ if the population is small. One can check that condition (\ref{ex:sweep})
can be replaced by a stronger one:
the measures $\mu_t$ given by $d\mu_t=P(t)f(x,i)\,m(dx,di)$ converge weakly to the measure
$\mu^*=p_0\delta_{(0,0)}+p_1\delta_{(0,1)}$.

Now, we assume additionally, that $g_i'(0)\ne 0$  for $i=0,1$ and $g_1'(a)\ne 0$
and  we check that stability and sweeping of
the semigroup  $\{P(t)\}_{t\ge 0}$ depends on the sign of the constant
\[
r_0=\frac{q_0(0)}{g'_0(0)}+\frac{q_1(0)}{g'_1(0)}.
\]
It is easy to check that both functions $\bar f_i$ are integrable in each interval outside the neighborhood  of
$0$. For any positive $\delta$ and  sufficiently small $x$ we have
\[
(r_0-\delta)x^{-1}\le r(x)\le (r_0+\delta)x^{-1}.
\]
It follows from these inequalities that there are some positive  numbers $c_1$, $c_2$ such that
\[
c_1 x^{-(r_0+\delta)} \le e^{-R(x)}\le c_2 x^{-(r_0-\delta)}
\]
for $x$ from a neighborhood  of $0$. Since $g_i(x)=g_i'(0)x+o(x)$, we obtain that  inequality (\ref{nier}) holds
when $r_0<0$ and that it  does not hold when $r_0>0$. Observe that $r_0<0$  iff
\begin{equation*}
q_0(0)g_1'(0)+q_1(0)g'_0(0)>0.
\end{equation*}
 This inequality  can be rewritten in the following way
\begin{equation}
\label{nier-int2}
\lambda =p_0g_0'(0)+p_1g'_0(0)>0.
\end{equation}
In the initial model we have $\lambda=b-\mu$, where $b=p_0b_0+p_1b_1$, and the number $\lambda$  can be
interpreted as the mean growth rate if the population is small. It explains why  the population becomes extinct
if $\lambda<0$ and  it survives if $\lambda>0$.
\end{example}


\section{Conclusions and summary}
\label{s:conclusion} In this paper we have presented a number of biological models described by PDMPs. The
models in Section \ref{sec:examples} have been chosen in such a way as to show that biological processes can
lead to various PDMPs, from a simple pure-jump Markov process with values in an Euclidean space to more advanced
Markov processes connected with individual based models in Section~\ref{ss:ppr-psm}. To study long-time
behaviour of  PDMPs we used the tool of stochastic semigroups on $L^1$-type spaces and their asymptotic
properties. We provided several examples of such semigroups in Section~\ref{s:stoch}. Theorems~\ref{asym-th2}
and~\ref{KTR} give criteria about asymptotic stability and sweeping with respect to compact sets of such
semigroups. Section~\ref{s:ltb} also contains examples of simple biological models which were used to illustrate
advanced techniques required to check that the related stochastic semigroup is asymptotically stable or
sweeping. Although these examples do not cover all models presented in Section~\ref{sec:examples}, the authors
believe that these results can be successfully applied to a wide range of models. In order to apply Theorem
\ref{KTR} one need to verify conditions (a) and (b), i.e., that the semigroup is irreducible and has some kernel
minorant. As we have mentioned in Section~\ref{s:ltb} one can check (b)  by using the H\"ormander's condition.
The final problem is to verify whether the semigroup is asymptotically stable if we already know that the
alternative  between  asymptotic stability and sweeping holds.  In more advanced models it might be very
difficult to prove  the existence of an invariant density in which case one can use the method of
Hasminski\u{\i} function (see \cite{RPT}) to exclude sweeping.

Our methods work quite well in the case of processes with switching dynamics or deterministic processes with
jumps if the jumps are "non-degenerated". An example of a "degenerated" jump is when we jump from a large part
of the phase space to one point. Such a "degenerated" jump appears in the neural activity model, when we jump
from points $(x,0,1)$, $x>\theta-a_E$, to the point $(0,0,2)$. Also in this model we have a "degenerated" jump
from the point $(0,t_R,2)$ to $(0,0,1)$ because $(0,0,1)$ is a stationary point of the related system of
differential equations and the process visits point $(0,0,1)$ with positive probability. But even in this case
one can induce a stochastic semigroup related to the PDMP if the measure $m$ on the phase space is an atom
measure at the point $(0,0,1)$ and the Lebesgue measure on the sets $\{(x,0,1) : x\in (-\infty,\theta)\}$
and  $\{(0,y,2) : y\in [0,t_R]\}$. We hope that it is possible to apply our technique to study the neural
activity model and to prove that the stochastic semigroup related to this model is asymptotically stable if
$a_E\lambda_E>a_I\lambda_I$. A priori our approach can not be applied to  processes connected with individual
based models from Section~\ref{ss:ppr-psm} where  it would be more convenient to work with more general
semigroups of probability measures. However, we are not aware of general results applicable in that example and
further work is required here.

\end{document}